\documentclass[12pt]{iopart}

\usepackage{iopams}

\let\ds=\displaystyle

\newcommand*{\hm}[1]{#1\nobreak\discretionary{}%
            {\hbox{$\mathsurround=0pt #1$}}{}}
\newtheorem{theorem}{Theorem}
\newtheorem{prop}{Proposition}
\newtheorem{lemma}{Lemma}
\newtheorem{remark}{Remark}

\newcommand{\End}{\mathop{\mathrm{End}}\nolimits}

\def\sk#1{\left(#1\right)}

\newcommand{\slt}{\mathfrak{sl}_2}
\newcommand{\glt}{\mathfrak{gl}_n}

\newcommand{\ra}{\rangle}
\newcommand{\La}{\Big\langle}
\newcommand{\Ra}{\Big\rangle}

\newcommand{\f}{f_\lambda}
\newcommand{\A}{{\mathcal A}}

\newcommand{\lfK}{{\cal K}}
\newcommand{\T}{\mathbb T}

\def\eqref#1{(\ref{#1})}
\def\ny{\nonumber}
\def\qed{\hfill$\Box$}
\def\sigmap{\mathcal{P}}

\begin{document}
\title[SOS model partition function and the elliptic weight functions]
{SOS model partition function and the elliptic\\ weight functions}
\author{S~Pakuliak$^{\dag\sharp}$, V Rubtsov$^{\ddag\sharp}$, A Silantyev$^{\dag\ddag}$}
\address{$^\dag$\ Laboratory of Theoretical Physics, JINR,
141980 Dubna, Moscow reg., Russia}
\address{$^\ddag$\ D\'epartment de Math\'ematiques, Universit\'e d'Angers,
2 Bd. Lavoisier, 49045 Angers, France}
\address{$^\sharp$\ Institute of Theoretical and Experimental Physics,
Moscow 117259, Russia}
\ead{ pakuliak@theor.jinr.ru, Volodya.Roubtsov@univ-angers.fr, silant@tonton.univ-angers.fr}
\bigskip

\begin{abstract}
We generalize a recent observation \cite{KhP} that the partition function of the 6-vertex
model with domain-wall boundary conditions can be obtained by
computing the projections of
the product of the total currents in the quantum affine algebra $U_{q}(\widehat{\mathfrak{sl}}_{2})$
in its current realization. A
generalization is proved for the the elliptic current algebra \cite{EF,ER1}.
The projections of the product of total currents are calculated explicitly and are
represented
as integral transforms of the product of the total currents. We
prove that the kernel of
this transform is proportional to the partition function of the SOS
model with
domain-wall boundary
conditions.
\end{abstract}
\hspace*{2.5cm}{\small \today}\hspace*{12mm}   \submitto{\JPA}
\vspace*{-11mm} 

\section{Introduction}

The main aim of this paper is to apply the method of elliptic current
projection to the computation of
the universal elliptic weight functions. The projection of currents first appeared in the works of
B. Enriquez and the second author~\cite{ER2}, \cite{ER3}, as
a method to construct a higher genus analog of the quantum groups in terms of Drinfeld
currents~\cite{D1}. The current (or ``new'') realization supplies a
quantum affine algebra with a second
co-product, the ``Drinfeld co-product''. The standard and Drinfeld co-products are related by a ``twist''
(see ~\cite{ER2}). The quantum algebra is decomposed in two different
ways a product of two Borel
 subalgebras.
For each subalgebra, we can consider its intersection with these two Borel subalgebras and
express it as their
product. Thus we obtain for each subalgebra a pair of projection operators from it to
 each of these intersections.
The above-mentioned twist is defined by a Hopf pairing of the subalgebras and the projection
operators. See Section 4 where
we recall an elliptic version of this construction.

S. Khoroshkin and the first author have applied this method
to a factorization of the universal $R$-matrix~\cite{DKhP} in
quantum affine algebras, in order to obtain
universal weight functions~\cite{KhP,EKhP} for arbitrary quantum affine algebras.
The weight functions play a fundamental role in the theory of deformed Knizhnik-Zamolodchikov and
Knizhnik-Zamolodchikov-Bernard equations.
In particular, in the case of $U_q(\widehat{\glt})$, acting by the projection of Drinfeld
currents onto the highest
weight vectors of irreducible finite-dimensional representations, one obtains exactly the
(trigonometric) weight functions
or off-shell Bethe vectors. In the canonical nested Bethe Ansatz, these objects are defined
implicitly by recursive relations. Calculations of the projections are an effective
way to determine the hierarchical relations of the nested Bethe Ansatz.

It was observed in \cite{KhP} that the projections for the
algebra $U_q(\widehat{\slt})$ can be represented as integral
transforms and that the kernels of these
transforms are proportional to the partition function of the finite 6-vertex model with
domain-wall boundary conditions (DWBC)~\cite{KhP}. We prove that the elliptic projections
described in~\cite{EF} make it possible to derive the partition function for elliptic models.
We show that the
calculation of the projections in the current elliptic algebra \cite{EF,ER1}
yields the partition function of   the Solid-On-Solid (SOS) model with domain-wall boundary
conditions.

The partition function for the finite 6-vertex model with domain wall boundary conditions was
obtained by Izergin~\cite{I87}, who derived recursion relations for
the partition function and solved them in determinant form.
The kernels of the projections satisfy the same recursion
relations and provide another formula for the partition function.

The problem of generalizing Izergin's determinant formula to the elliptic case
has been extensively discussed in the last two decades.
One can prove that the statistical sum of the SOS model with DWBC  cannot be represented in
the form of a single determinant. While this paper was in preparation,
H. Rosengren \cite{Ros08} showed
that this statistical sum for an $n\times n$ lattice can be written as
a sum of $2^n$ determinants, thus
generalizing Izergin's determinant formula. His approach relates to some dynamical generalization
of the method of
Alternating-Sign Matrices and follows the famous
combinatorial proof of Kuperberg \cite{Kup}.

We expect that the projection method gives a universal form for the
elliptic weight function \cite{TV} as it does for the quantum affine
algebras \cite{KPT}. When this universal weight function is represented as
an integral transform of the product of the elliptic currents, we show that the kernel
of this transform gives an expression of the partition function for the SOS model.
On the one hand we generalize Izergin's recurrent relations and on the
other hand we generalize to the elliptic case the method
proposed in \cite{KhP} for calculating the projections.
We check that the kernel extracted from the universal weight function
and multiplied by a certain factor satisfies
the recursion relations that have been obtained,
which uniquely define the partition function for the SOS model with DWBC.
Our formula given by the projection method coincides with Rosengren's.

An interesting open problem which deserves more a extensive study is
the relation of the projection method
with the elliptic Sklyanin-Odesskii-Feigin algebras. It was observed in the pioneering paper \cite{ER1} that
half of the elliptic current generators satisfy the commutation
relations of the $W$-elliptic algebras of
Feigin. Another intriguing relation was observed in \cite{FO}: there existe a certain subalgebra in the
``$\lambda$-generalization''
of the Sklyanin algebra
such that its generators obey the Felder's $R$-matrix
quadratic relations given in \cite{FVT}. The latter paper gives also a description of the elliptic
Bethe eigenvectors (the elliptic weight functions).

This is a strong indication that the projection method should be considered and interpreted in the
framework of the (generalized) Sklyanin-Odesskii-Feigin algebras. We
hope to discuss this problem elsewhere.

The main results of this paper were reported at the 7-th International
Workshop on
``Supersymmetry and Quantum
Symmetry'' in JINR, Dubna (Russia), July 30 - August 4, 2007.

The paper is organized as follows. In section~\ref{sec2} we briefly review the finite 6-vertex model
with DWBC, and we present the formulae for the partition function: Izergin's determinant formula
and the formula obtained by the projection method. Section~\ref{sec3} is devoted to the SOS model
with DWBC. We briefly introduce the model and pose the problem of how to calculate
the partition function of this model. We derive analytical properties of the partition function
that allow us to reconstruct the partition function exactly.
In section~\ref{sec4} we introduce the projections in terms of the currents
for the elliptic algebra, following~\cite{EF}. We generalize the method proposed in~\cite{KhP}
to this case in order to obtain the integral representation of the projections of products of currents.
Then, using a Hopf pairing, we extract the kernel and show
that it  satisfies all the necessary analytical properties
of the partition function of the SOS model with DWBC.
In Section~\ref{sec5}, we investigate the trigonometric degeneration of the elliptic model
and of the partition function with DWBC. We arrive at the 6-vertex model case in two steps.
The model obtained after the first step is a trigonometric SOS model. Then we show that
the degeneration of the
expression derived in Section~\ref{sec4} coincides with the known expression for the
6-vertex model partition function with DWBC.
An appendix contains the necessary information on the properties of elliptic polynomials.

\section{Partition function of the finite 6-vertex model}
\label{sec2}

\begin{figure}[h]
\begin{center}
\begin{picture}(180,130)
\put(05,63){$a(z,w)$}
\put(05,-15){$a(z,w)$}
\put(23,118){$+$}
\put(0,104){$+$}
\put(33,104){$+$}
\put(23,78){$+$}
\put(40,100){\vector(-1,0){40}}
\put(40,100){\vector(-1,0){3}}
\put(20,120){\vector(0,-1){40}}
\put(20,120){\vector(0,-1){3}}

\put(85,63){$b(z,w)$}
\put(85,-15){$b(z,w)$}
\put(103,118){$+$}
\put(113,104){$-$}
\put(80,104){$-$}
\put(103,78){$+$}
\put(80,100){\vector(1,0){40}}
\put(80,100){\vector(1,0){3}}
\put(100,120){\vector(0,-1){40}}
\put(100,120){\vector(0,-1){3}}

\put(165,63){$c(z,w)$}
\put(165,-15){$\bar c(z,w)$}
\put(183,118){$-$}
\put(193,104){$+$}
\put(160,104){$-$}
\put(183,78){$+$}
\put(160,100){\line(1,0){40}}
\put(200,100){\vector(-1,0){3}}
\put(160,100){\vector(1,0){3}}
\put(180,80){\vector(0,1){40}}
\put(180,80){\vector(0,-1){3}}

\put(23,38){$-$}
\put(0,24){$-$}
\put(33,24){$-$}
\put(23,00){$-$}
\put(0,20){\vector(1,0){40}}
\put(0,20){\vector(1,0){3}}
\put(20,00){\vector(0,1){40}}
\put(20,00){\vector(0,1){3}}

\put(103,38){$-$}
\put(113,24){$+$}
\put(80,24){$+$}
\put(103,00){$-$}
\put(120,20){\vector(-1,0){40}}
\put(120,20){\vector(-1,0){3}}
\put(100,00){\vector(0,1){40}}
\put(100,00){\vector(0,1){3}}

\put(183,38){$+$}
\put(193,24){$-$}
\put(160,24){$+$}
\put(183,00){$-$}
\put(200,20){\line(-1,0){40}}
\put(200,20){\vector(1,0){3}}
\put(160,20){\vector(-1,0){3}}
\put(180,00){\line(0,1){40}}
\put(180,00){\vector(0,1){3}}
\put(180,40){\vector(0,-1){3}}
\end{picture}
\medskip
\end{center}
\caption{\footnotesize  Graphical presentation of the Boltzmann weights.}
\label{fig1}
\end{figure}
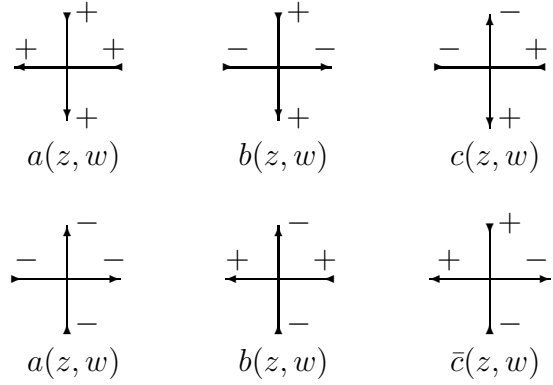

Let us consider a statistical system on a square $n\times n$ lattice, where the
columns and rows are numbered from $1$ to $n$ from right to left and
from bottom to top, respectively.
This is a 6-vertex model where the vertices on the lattice are associated with Boltzmann weights
which depend on the configuration of the arrows around a given vertex.
The six possible configurations are shown in Fig.~\ref{fig1}, The weights are functions of
two spectral parameters $z$, $w$ and ananisotropy parameter $q$:
\begin{equation} \label{Bw}
\eqalign{a(z,w)&=qz-q^{-1}w, \qquad  b(z,w)=z-w, \cr
c(z,w)&=(q-q^{-1})z, \qquad \bar c(z,w)=(q-q^{-1})w.}
\end{equation}

Let us associate the sign `$+$' to the arrows directed upward  and to
the left,
while the  sign `$-$' is associated to to the arrows directed downward  and to
the right as
shown in Fig.~\ref{fig1}. The Boltzmann weights~\eqref{Bw} are gathered in the matrix
\begin{equation}\label{Bw1}
R(z,w)=\sk{\begin{array}{cccc}a(z,w)&0&0&0\\ 0&b(z,w)& \bar c(z,w)&0\\
0&c(z,w)&b(z,w)&0\\0&0&0&a(z,w)
\end{array}}
\end{equation}
acting in the space $\mathbb C^2\otimes\mathbb C^2$ with the basis $e_{\alpha}\otimes e_{\beta}$,
$\alpha,\beta=\pm$. The entry $R(z,w)^{\alpha\beta}_{\gamma\delta}$,
$\alpha,\beta,\gamma,\delta=\pm$ coincides with the Boltzmann weight corresponding to
Fig.~\ref{fig2}:

\begin{figure}[h]
\begin{center}
\begin{picture}(40,40)
\put(23,28){$\alpha$}
\put(0,14){$\delta$}
\put(33,14){$\beta$}
\put(23,-10){$\gamma$}
\put(40,10){\line(-1,0){40}}
\put(40,10){\line(-1,0){3}}
\put(20,-10){\line(0,1){40}}
\put(20,-10){\line(0,1){3}}
\end{picture}
\end{center}
\caption{\footnotesize  The Boltzmann weight $R(z,w)^{\alpha\beta}_{\gamma\delta}$.}
\label{fig2}
\end{figure}
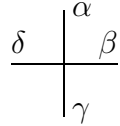

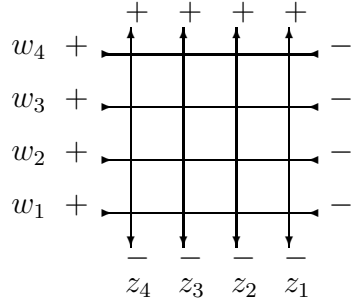
\begin{figure}
\begin{center}
\begin{picture}(180,100)
\put(50,70){\line(1,0){80}}
\put(50,70){\vector(1,0){3}}
\put(130,70){\vector(-1,0){3}}
\put(15,70){$w_4$} \put(35,70){$+$} \put(135,70){$-$}
\put(50,50){\line(1,0){80}}
\put(50,50){\vector(1,0){3}}
\put(130,50){\vector(-1,0){3}}
\put(15,50){$w_3$} \put(35,50){$+$} \put(135,50){$-$}
\put(50,30){\line(1,0){80}}
\put(50,30){\vector(1,0){3}}
\put(130,30){\vector(-1,0){3}}
\put(15,30){$w_2$} \put(35,30){$+$} \put(135,30){$-$}
\put(50,10){\line(1,0){80}}
\put(50,10){\vector(1,0){3}}
\put(130,10){\vector(-1,0){3}}
\put(15,10){$w_1$} \put(35,10){$+$} \put(135,10){$-$}

\put(60,0){\vector(0,1){80}}
\put(60,0){\vector(0,-1){3}}
\put(58,-10){$-$} \put(58,-20){$z_4$} \put(58,83){$+$}
\put(80,0){\vector(0,1){80}}
\put(80,0){\vector(0,-1){3}}
\put(78,-10){$-$} \put(78,-20){$z_3$} \put(78,83){$+$}
\put(100,0){\vector(0,1){80}}
\put(100,0){\vector(0,-1){3}}
\put(98,-10){$-$} \put(98,-20){$z_2$} \put(98,83){$+$}
\put(120,0){\vector(0,1){80}}
\put(120,0){\vector(0,-1){3}}
\put(118,-10){$-$} \put(118,-20){$z_1$} \put(118,83){$+$}
\end{picture}
\medskip
\end{center}
\caption{\footnotesize Inhomogeneous lattice with domain wall boundary conditions.}
\label{fig3}
\end{figure}

Different repartitions of the arrows on the edges form different configurations $\{C\}$.
A Boltzmann weight of the lattice is the product of the Boltzmann weights at each vertex.
We define the partition function of the model as the sum of
the Boltzmann weights of the lattice over all possible configurations, subject
to some boundary conditions:
\begin{equation}\label{sts1}
Z(\{z\},\{w\})=\sum_{\{C\}}\prod_{i,j=1}^n R(z_i,w_j)^{\alpha_{ij}\beta_{ij}}_{\gamma_{ij}\delta_{ij}}.
\end{equation}
Here $\alpha_{ij}$, $\beta_{ij}$, $\gamma_{ij}$, $\delta_{ij}$ are
the signs corresponding to the arrows
around the $(i,j)$-th vertex. We consider an inhomogeneous model where the Boltzmann
weights depend on the column by the variable $z_i$ and on the row by the variable $w_j$
(see Fig.~\ref{fig3}).

We choose the so-called domain-wall boundary conditions (DWBC) that fix the boundary arrows (signs)
as shown in Fig.~\ref{fig3}. In other words, the arrows enter on the left and
right boundaries and leave on the lower and upper ones.

In \cite{I87}, A.\,G.\,Izergin found a determinant representation for the partition
function of the lattice with DWBC,
\begin{eqnarray}
\fl Z(\{z\},\{w\})= (q-q^{-1})^n \prod_{m=1}^n w_m \times  \ny\\
\fl\qquad \times\frac{\ds \prod_{i,j=1}^n(z_i-w_j)(qz_i-q^{-1}w_j)}
{\ds \prod_{n\ge i>j\ge1}(z_i-z_j)(w_j-w_i)}
\det\left|\left|\frac{1}{(z_i-w_j)(qz_i-q^{-1}w_j)}\right|\right|_{i,j=1,\ldots,n}\, .\label{stat-s}
\end{eqnarray}
Izergin's idea was to prove a symmetry of the polynomial~\eqref{sts1}, and then use it
to find recursion relations for the quantity $Z(\{z\},\{w\})$ and to observe that these
recursion relations allow the reconstruction of  $Z(\{z\},\{w\})$ in a
unique way and that the
same recursion relations are valid for the determinant formula~\eqref{stat-s}.

On the other hand it was observed that the kernel of the projection of $n$ currents is a
polynomial of the same degree, and satisfies the same recursion relations~\cite{KhP}.
It means that this kernel coincides with the partition function for
the $n\times n$ lattice. Moreover,
the theory of projections gives another expression for the partition function:
\begin{eqnarray}
\fl Z(\{z\},\{w\}) = (q-q^{-1})^n \prod_{m=1}^n w_m \prod_{n\ge i>j\ge1}\frac{q^{-1}w_i-qw_j}
 {w_i-w_j}\times \ny \\
\fl\qquad \times \sum_{\sigma\in S_n}\prod_{1\le i<j\le n \atop \sigma(i)>\sigma(j)}
  \frac{qw_{\sigma(i)}-q^{-1}w_{\sigma(j)}}{q^{-1}w_{\sigma(i)}-qw_{\sigma(j)}}
  \prod_{n\ge i>k\ge1}(qz_i-q^{-1}w_{\sigma(k})\prod_{1\le i<k\le n} (z_i-w_{\sigma(k)}),
  \label{stat-s_pr}
\end{eqnarray}
where $S_n$ is the group of permutations. Here the factor
$\frac{qw_{\sigma(i)}-q^{-1}w_{\sigma(j)}}{q^{-1}w_{\sigma(i)}-qw_{\sigma(j)}}$
appears in the product if both conditions $i<j$ and $\sigma(i)>\sigma(j)$ are
satisfied simultaneously.

\section{Partition function for the SOS model}
\label{sec3}

\subsection{Description of the SOS model}
\label{sec31}

The SOS model is a face model. We introduce it in terms of heights as
usual, but then we represent
it in the $R$-matrix formalism as in~\cite{FS}. This language is more convenient to
generalize the results reviewed in
Section~\ref{sec2} and to prove the symmetry of the partition function.

\begin{figure}
\begin{center}
\begin{picture}(150,90)
\put(50,70){\line(1,0){80}}
\put(35,76){$n$}
\put(50,50){\line(1,0){80}}
\put(35,60){$\ldots$}
\put(50,30){\line(1,0){80}}
\put(5,38){$j=$}
\put(35,38){$2$}
\put(50,10){\line(1,0){80}}
\put(35,19){$1$}
\put(35,00){$0$}
\put(60,0){\line(0,1){80}}
\put(48,-15){$n$}
\put(80,0){\line(0,1){80}}
\put(68,-15){$\ldots$}
\put(100,0){\line(0,1){80}}
\put(88,-15){$2$}
\put(120,0){\line(0,1){80}}
\put(108,-15){$1$}
\put(123,-15){$0$}
\put(20,-15){$i=$}
\end{picture}
\end{center}
\caption{\footnotesize The numbering of faces.}
\label{fig5}
\end{figure}
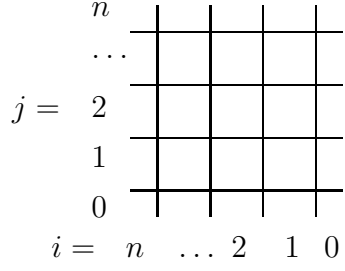

Consider a square $n\times n$ lattice with vertices enumerated by the index $i=1,\ldots,n$ as
in the previous case. It has $(n+1)\times(n+1)$ faces enumerated by pairs $(i,j)$, $i,j=0,\ldots,n$
(see Fig.~\ref{fig5}). To each face we assign a complex number called
its height in such a way that the
differences of the heights corresponding to the neighboring faces are $\pm1$. Let us denote
by $d_{ij}$ the height corresponding to the face $(i,j)$ placed to the
upper left of the
vertex $(i,j)$. Then the last condition can be written in the form $|d_{ij}-d_{i-1,j}|=1$,
for $i=1,\ldots,n$, $j=0,\ldots,n$, and $|d_{ij}-d_{ij-1}|=1$ for $i=0,\ldots,n$, $j=1,\ldots,n$.
Each distribution of heights $d_{ij}$ ($i,j=0,\ldots,n$) subject to
these conditions and to
boundary conditions defines a configuration of the model. It means that the partition
function of this model can be represented in the form
\begin{equation}
 Z=\sum_{C}\prod_{i,j=1}^n W_{ij}(d_{i,j-1},d_{i-1,j-1},d_{i-1,j},d_{ij}), \label{pf_def}
\end{equation}
where $W_{ij}(d_{i,j-1},d_{i-1,j-1},d_{i-1,j},d_{ij})$ is the
Boltzmann weight of the $(i,j)$-th vertex
depending on the configuration by means of the heights of the
neighboring faces
%??
as follows~\cite{J}
\begin{equation}\label{BWa}
\eqalign{W_{ij}(d+1,d+2,d+1,d)=a(u_i-v_j)=\theta(u_i-v_j+\hbar),  \cr
 W_{ij}(d-1,d-2,d-1,d)=a(u_i-v_j)=\theta(u_i-v_j+\hbar), \cr
 W_{ij}(d-1,d,d+1,d)=b(u_i-v_j;\hbar d)=\frac{\theta(u_i-v_j)\theta(\hbar d+\hbar)}{\theta(\hbar d)}, \cr
 W_{ij}(d+1,d,d-1,d)=\bar b(u_i-v_j;\hbar d)=\frac{\theta(u_i-v_j)\theta(\hbar d-\hbar)}{\theta(\hbar d)}, \cr
 W_{ij}(d-1,d,d-1,d)=c(u_i-v_j;\hbar d)=\frac{\theta(u_i-v_j+\hbar d)\theta(\hbar)}{\theta(\hbar d)}, \cr
 W_{ij}(d+1,d,d+1,d)=\bar c(u_i-v_j;\hbar d)=\frac{\theta(u_i-v_j-\hbar d)\theta(\hbar)}{\theta(-\hbar d)}.}
\end{equation}
As in the 6-vertex case the variables $u_i$, $v_j$ are attached to the $i$-th vertical and $j$-th horizontal
lines respectively, $\hbar$ is a nonzero anisotropy parameter~\footnote{In the elliptic case, we
use additive variables $u_i$, $v_j$ and an additive anisotropy parameter $\hbar$ instead of the
multiplicative variables $z_i=e^{2\pi iu_i}$, $w_i=e^{2\pi iv_i}$ and
the multiplicative parameter
$q=e^{\pi i\hbar}$.}. The weights are expressed by means of
the ordinary odd theta-function defined by the conditions
\begin{equation}
 \theta(u+1)=-\theta(u),\quad  \theta(u+\tau)=-e^{-2\pi iu-\pi i\tau}\theta(u),\quad  \theta'(0)=1.
 \label{theta}
\end{equation}

\begin{figure}
\begin{center}
 \begin{picture}(240,170)
% left-down
\put(00,30){\line(1,0){60}}\put(30,00){\line(0,1){60}}
\put(-10,30){\bfseries $-$}\put(60,30){$-$}\put(30,-10){\bfseries $-$}\put(30,60){\bfseries $-$}
\put(0,10){\small $d-1$}\put(10,40){\small $d$}\put(35,40){\small $d-1$}\put(35,10){\small $d-2$}
\put(10,-20){\small  $a(u_i-v_j)$}
% left-up
\put(00,130){\line(1,0){60}}\put(30,100){\line(0,1){60}}
\put(-10,130){\bfseries $+$}\put(60,130){\bfseries $+$}\put(30,90){\bfseries $+$}\put(30,160){\bfseries $+$}
\put(0,110){\small $d+1$}\put(10,140){\small $d$}\put(35,140){\small $d+1$}\put(35,110){\small $d+2$}
\put(10,77){\small  $a(u_i-v_j)$}
% centre-down
\put(100,30){\line(1,0){60}}\put(130,00){\line(0,1){60}}
\put(90,30){\bfseries $+$}\put(160,30){$+$}\put(130,-10){\bfseries $-$}\put(130,60){\bfseries $-$}
\put(100,10){\small $d+1$}\put(110,40){\small $d$}\put(135,40){\small $d-1$}\put(145,10){\small $d$}
\put(105,-20){\small  $\bar b(u_i-v_j;d)$}
% centre-up
\put(100,130){\line(1,0){60}}\put(130,100){\line(0,1){60}}
\put(90,130){\bfseries $-$}\put(160,130){\bfseries $-$}\put(130,90){\bfseries $+$}\put(130,160){\bfseries $+$}
\put(100,110){\small $d-1$}\put(110,140){\small $d$}\put(135,140){\small $d+1$}\put(145,110){\small $d$}
\put(105,77){\small  $b(u_i-v_j;d)$}
% right-down
\put(200,30){\line(1,0){60}}\put(230,00){\line(0,1){60}}
\put(190,30){\bfseries $+$}\put(260,30){$-$}\put(230,-10){\bfseries $-$}\put(230,60){\bfseries $+$}
\put(200,10){\small $d+1$}\put(210,40){\small $d$}\put(235,40){\small $d+1$}\put(245,10){\small $d$}
\put(205,-20){\small  $\bar c(u_i-v_j;d)$}
% right-up
\put(200,130){\line(1,0){60}}\put(230,100){\line(0,1){60}}
\put(190,130){\bfseries $-$}\put(260,130){\bfseries $+$}\put(230,90){\bfseries $+$}\put(230,160){\bfseries $-$}
\put(200,110){\small $d-1$}\put(210,140){\small $d$}\put(235,140){\small $d-1$}\put(245,110){\small $d$}
\put(205,77){\small  $c(u_i-v_j;d)$}
  \end{picture}
  \bigskip
 \end{center}
\caption{\footnotesize The Boltzmann weights for the SOS model.}
\label{fig4}
\end{figure}

Let us introduce the notations
\begin{equation}
\fl \alpha_{ij}=d_{i-1,j}-d_{ij},\  \beta_{ij}=d_{i-1,j-1}-d_{i-1,j},\
  \gamma_{ij}=d_{i-1,j-1}-d_{i,j-1},\  \delta_{ij}=d_{i,j-1}-d_{ij}. \label{albegade}
\end{equation}
The differences~\eqref{albegade} take the values $\pm1$ and we attach them to the corresponding edges
as in Fig.~\ref{fig2}: $\gamma_{i,j+1}=\alpha_{ij}$ is the sign attached to the vertical edge
connecting the $(i,j)$-th vertex to the $(i,j+1)$-st one, $\beta_{i,j+1}=\delta_{ij}$ is the sign
attached to the horizontal edge connecting the $(i,j)$-th vertex to the $(i+1,j)$-th one.
The configuration can be considered as a distribution of these signs on the internal edges
subject to the conditions $\alpha_{ij}+\beta_{ij}=\gamma_{ij}+\delta_{ij}$, $i,j=1,\ldots,n$.
In terms of signs on the external edges the DWBC are the same as shown in Fig.~\ref{fig3}.
Additionally, we have to fix one of the boundary heights, for example, $d_{nn}$.

The Boltzmann weights~\eqref{BWa} can be represented as the entries of
a dynamical elliptic
$R$-matrix~\cite{FS}:
\begin{eqnarray}
 W_{ij}(d_{i,j-1},d_{i-1,j-1},d_{i-1,j},d_{ij})
   =R(u_i-v_j;\hbar d_{ij})^{\alpha_{ij}\beta_{ij}}_{\gamma_{ij}\delta_{ij}}, \ny \\
 R(u;\lambda)=\left(\begin{array}{cccc}
   a(u) & 0              & 0              & 0              \\
  0              & b(u;\lambda) & \bar c(u;\lambda)      & 0              \\
  0              & c(u;\lambda) & \bar b(u;\lambda) & 0              \\
  0              & 0              & 0              & a(u)
 \end{array}\right). \label{Rz}
\end{eqnarray}

Let $\T^{\alpha_{in}}_{\gamma_{i1}}(u_i,\{v\},\lambda_i)$ be the {\itshape column transfer matris}.
It is a  matrix-valued
function of $u_i$, all spectral parameters
$v_j$, $j=1,\ldots,n$ and the parameters $\lambda_i$ related to the heights:
\begin{eqnarray}
\fl \T^{\alpha_{in}}_{\gamma_{i1}}(u_i,\{v\},\lambda_i)^{\beta_{in}\ldots\beta_{i1}}_{\delta_{in}\ldots
 \delta_{i1}} = \label{colTM} \\
\fl = \Big(R^{(n+1,n)}(u_i-v_n;\lambda_{in})R^{(n+1,n-1)}(u_i-v_{n-1};\lambda_{i,n-1})\cdots
    R^{(n+1,1)}(u_i-v_1;\lambda_{i1})\Big)^{\alpha_{in};\beta_{in}\ldots\beta_{i1}}_{\gamma_{i1};
    \delta_{in}\ldots\delta_{i1}} \ny \\
\fl =  \Big(R^{(n+1,n)}(u_i-v_n;\Lambda_{in})R^{(n+1,n-1)}(u_i-v_{n-1};\Lambda_{i,n-1})\cdots
    R^{(n+1,1)}(u_i-v_1;\Lambda_{i1})\Big)^{\alpha_{in};\beta_{in}\ldots\beta_{i1}}_{\gamma_{i1};
    \delta_{in}\ldots\delta_{i1}}, \ny
\end{eqnarray}
where $\lambda_{ij}=\hbar d_{ij}=\lambda_i+\hbar\sum\limits_{l=j+1}^n\delta_{il}$,
$\lambda_i=\hbar d_{in}=\lambda+\hbar\sum\limits_{l=i+1}^n\alpha_{ln}$,
$\Lambda_{ij}=\lambda_i+\hbar\sum\limits_{l=j+1}^n H^{(l)}$.
The matrix $H^{(l)}$ acts in the $l$-th two-dimensional space  $V_{l}\cong\mathbb C^2$
as   ${\rm diag}(1,-1)$  and  the $R$-matrix $R^{(a,b)}$ acts nontrivially in the
         tensor product $V_a\otimes V_b$.
The superscript $n+1$ in the $R$-matrices  is regarded as belonging
%?
to an auxiliary space
         $V_{n+1}\cong\mathbb C^2$. The partition function~\eqref{pf_def} corresponding
          to DWBC ($\alpha_{in}=+1$,  $\beta_{1i}=-1$, $\gamma_{i1}=-1$, $\delta_{ni}=+1$,
$i=1,\ldots,n$) can be represented by means of the column transfer matrices:
\begin{equation}
 Z^{+-}_{-+}(\{u\},\{v\};\lambda)
=\Big(\T^+_-(u_1,\{v\},\lambda_1)\cdots\T^+_-(u_n,\{v\},\lambda_n)\Big)^{-\ldots-}_{+\ldots+}, \label{ZcT}
\end{equation}
where $\lambda_i=\lambda+\hbar(n-i)$. Similarly one can define the row transfer matrix.

\subsection{Analytical properties of the partition function}
\label{sec32}

Here we describe the analytical properties of the SOS model partition
function which are analogous to those used by
A.G. Izergin in order to recover the partition function of  the 6-vertex model. These properties
uniquely define this partition function.

\begin{prop} \label{propZsym}
 The partition function with DWBC $Z^{+-}_{-+}(\{u\},\{v\};\lambda)$ is a symmetric function in
 both sets of variables $u_i$ and $v_j$.
\end{prop}

\noindent
The proof is based on the dynamical Yang-Baxter equation (DYBE) for the $R$-matrix~\cite{FS}
\begin{eqnarray}
 R^{(12)}(t_1-t_2;\lambda)R^{(13)}(t_1-t_3;\lambda+\hbar H^{(2)})R^{(23)}(t_2-t_3;\lambda)= \ny\\
=R^{(23)}(t_2-t_3;\lambda+\hbar H^{(1)})R^{(13)}(t_1-t_3;\lambda)R^{(12)}(t_1-t_2;\lambda+\hbar H^{(3)}).\ny
\end{eqnarray}
In order to prove the symmetry of the partition function $Z^{+-}_{-+}(\{u\},\{v\};\lambda)$ under
the permutation $v_j\leftrightarrow v_{j-1}$, we rewrite the DYBE in the form
\begin{eqnarray}
 &R^{(n+1,j)}(u_i-v_j;\Lambda_{ij})R^{(n+1,j-1)}(u_i-v_{j-1};\Lambda_{ij}+\hbar H^{(j)})\ny \\
 &\quad \times R^{(j,j-1)}(v_j-v_{j-1};\Lambda_{ij})= R^{(j,j-1)}(v_j-v_{j-1};\Lambda_{ij}+\hbar H^{(n+1)})
\label{DYBEuse} \\
 &\quad \times R^{(n+1,j-1)}(u_i-v_{j-1};\Lambda_{ij})
 R^{(n+1,j)}(u_i-v_j;\Lambda_{ij}+\hbar H^{(j-1)}). \ny
\end{eqnarray}
 Multiplying the $i$-th
column matrix~\eqref{colTM} by $R^{(j,j-1)}(v_j-v_{j-1};\Lambda_{ij})$
to the right and moving it
to the left using~\eqref{DYBEuse}, the relation $[H_1+H_2,R(u,\lambda)]=0$ and the equality
$\Lambda_{ij}+\hbar\alpha_{in}=\Lambda_{i-1,j}$, we obtain
\begin{eqnarray}
\fl \T^{\alpha_{in}}_{\gamma_{i1}}(u_i,\{v\},\lambda_i)R^{(j,j-1)}(v_j-v_{j-1};\Lambda_{ij})
 =R^{(j,j-1)}(v_j-v_{j-1};\Lambda_{i-1,j})\ny\\
\fl\quad \times \sigmap^{(j,j-1)}
 \T^{\alpha_{in}}_{\gamma_{i1}}(u_i,\{v_j\leftrightarrow v_{j-1}\},\lambda_i)
 \cdots \T^{\alpha_{nn}}_{\gamma_{n1}}(u_n,\{v_j\leftrightarrow v_{j-1}\},\lambda_n)
 \sigmap^{(j,j-1)}, \label{TR_RT}
\end{eqnarray}
where $\sigmap\in\End(\mathbb C^2\otimes\mathbb C^2)$ is a permutation matrix:
$\sigmap(e_1\otimes e_2)=e_2\otimes e_1$ for all $e_1,e_2\in\mathbb
C^2$ and the
notation $\{v_j\leftrightarrow v_{j-1}\}$ means that the set of parameters $\{v\}$ with
$v_{j-1}$ and $v_j$ are interchanged. Multiplying the product
of the column matrix by $R^{(j,j-1)}(v_j-v_{j-1};\Lambda_{nj})$ to the right and moving it to the
left using~\eqref{TR_RT} one yields
\begin{eqnarray}
\fl \T^{\alpha_{1n}}_{\gamma_{11}}(u_1,\{v\},\lambda_1)\cdots\T^{\alpha_{nn}}_{\gamma_{n1}}(u_n,\{v\},
 \lambda_n)R^{(j,j-1)}(v_j-v_{j-1};\Lambda_{nj}) =R^{(j,j-1)}(v_j-v_{j-1};\Lambda_{0,j})\ny \\
\fl \quad\times
  \sigmap^{(j,j-1)}
 \T^{\alpha_{1n}}_{\gamma_{11}}(u_1,\{v_j\leftrightarrow v_{j-1}\},\lambda_1)\cdots
 \T^{\alpha_{nn}}_{\gamma_{n1}}(u_n,\{v_j\leftrightarrow v_{j-1}\},\lambda_n)\sigmap^{(j,j-1)},  \label{RTTT}
\end{eqnarray}
where $\Lambda_{0,j}=\lambda+\hbar\sum\limits_{i=1}^n\alpha_{in}+\hbar
\sum\limits_{l=j+1}^n H^{(l)}$, $\Lambda_{nj}=\lambda_n=\lambda$. Finally,
comparing the matrix element $\big(\cdot\big)^{-,\ldots,-}_{+,\ldots,+}$ of both
sides of~\eqref{RTTT}, taking into account Formula~\eqref{ZcT} and the identities
\begin{equation*}
\fl  R(u,\lambda)^{--}_{\gamma\delta}=a(u)\delta^-_\gamma\delta^-_\delta,
\quad R(u,\lambda)^{\alpha\beta}_{++}=a(u)\delta^\alpha_+\delta^\beta_+,
 \quad  \sigmap^{--}_{\gamma\delta}=\delta^-_\gamma\delta^-_\delta,
  \quad \sigmap^{\alpha\beta}_{++}=\delta^\alpha_+\delta^\beta_+,
\end{equation*}
(where $\delta^\alpha_\gamma$ is the Kronecker symbol) and substituting $\alpha_{in}=+1$,
$\gamma_{i1}=-1$, one derives
\begin{equation}
 Z^{+-}_{-+}(\{u\},\{v\};\lambda)=Z^{+-}_{-+}(\{u\},\{v_j\leftrightarrow v_{j-1}\};\lambda). \label{Z_perm_w}
\end{equation}
Similarly, using the row transfer matrix one can obtain the following
equality from the DYBE:
\begin{equation}
  Z^{+-}_{-+}(\{u\},\{v\};\lambda)=Z^{+-}_{-+}(\{u_j\leftrightarrow u_{j-1}\},\{v\};\lambda). \label{Z_perm_z}
\end{equation}
The partition function with DWBC satisfies relations~\eqref{Z_perm_w}, \eqref{Z_perm_z}
for each $j=1,\ldots,n$, which is sufficient to establish the symmetry under an arbitrary permutation. \qed

\begin{prop} \label{propZep}
 The partition function with DWBC~\eqref{ZcT} is an elliptic polynomial
 \footnote{The definition of elliptic polynomials and their properties are
given in \ref{Appendix_ep}.} of degree $n$
 with character $\chi$ in each variable $u_i$, where
\begin{equation}
 \chi(1)=(-1)^n,\qquad
 \chi(\tau)=(-1)^n\exp\Big(2\pi i\big(\lambda+\sum_{j=1}^n v_j\big)\Big). \label{char_chi}
\end{equation}
\end{prop}

\noindent
Due to the symmetry
%of what ?
with respect to the variables $\{u\}$ it is
sufficient to prove the proposition for the variable $u_n$. To represent explicitly the
dependence of $Z^{+-}_{-+}(\{u\},\{v\};\lambda)$ on $u_n$, we consider
all the possibilities for the
states of the
edges attached to the vertices located in the $n$-th column. First, consider the $(n,n)$-th vertex.
Due to the boundary conditions $\alpha_{nn}=\delta_{nn}=+1$ and to the condition
$\alpha_{nn}+\beta_{nn}=\gamma_{nn}+\delta_{nn}$ we have two possibilities: either
$\beta_{nn}=\gamma_{nn}=-1$ or $\beta_{nn}=\gamma_{nn}=+1$. In the first case, one has a
unique possibility for the rest of the $n$-th column: $\gamma_{nj}=-1$,
$\beta_{nj}=+1$, $j=1,\ldots,n-1$; in the second case, there are two possibilities for the $(n,n-1)$-st
vertex: either $\beta_{n,n-1}=\gamma_{n,n-1}=-1$ or $\beta_{n,n-1}=\gamma_{n,n-1}=+1$, etc.
Finally the partition function is represented in the from
\begin{eqnarray}
 Z^{+-}_{-+}(\{u\},\{v\};\lambda)=
  \sum_{k=1}^n \prod_{j=k+1}^n  a(u_n-v_j)\;\bar c(u_n-v_k;\lambda+(n-k)\hbar) \ny \\
\quad  \times\prod_{j=1}^{k-1} \bar b(u_n-v_j;\lambda+(n-j)\hbar) g_k(u_{n-1},\ldots,u_1,\{v\};\lambda), \ny
\end{eqnarray}
where $g_k(u_{n-1},\ldots,u_1,\{v\};\lambda)$ are functions which do
not depend on $u_n$.
Each term of this sum is an elliptic polynomial of degree $n$ with the same
character~\eqref{char_chi} in the variable $u_n$. \qed
\medskip

%{\bfseries Remark \remnumber\label{rem1}.}
\begin{remark}\label{rem1} Similarly, one can prove that the function
$Z^{+-}_{-+}(\{u\},\{v\};\lambda)$ is an elliptic polynomial of degree $n$ with
character $\tilde\chi$ in each variable $v_i$, where $\tilde\chi(1)=(-1)^n$, \\
$\tilde\chi(\tau)=(-1)^n e^{2\pi i(-\lambda+\sum_{i=1}^n u_i)}$.
\end{remark}

\begin{prop} \label{propZrec}
The $n$-th partition function with DWBC~\eqref{ZcT} with the condition
$u_n=v_n-\hbar$ can be expressed through the $(n-1)$-st partition function:
\begin{eqnarray}
\fl Z^{+-}_{-+}(u_n=v_n-\hbar,u_{n-1},\ldots,u_1;v_n,v_{n-1},\ldots,v_1;\lambda)= \label{Zrec} \\
\fl =\frac{\theta(\lambda+n\hbar)\theta(\hbar)}{\theta(\lambda+(n-1)\hbar)}
  \prod_{m=1}^{n-1}\Big(\theta(v_n-v_m-\hbar)\theta(u_m-v_n)\Big)\,
 Z^{+-}_{-+}(u_{n-1},\ldots,u_1;v_{n-1},\ldots,v_1;\lambda). \ny
\end{eqnarray}
\end{prop}

\noindent
Considering the $n$-th column and the $n$-th row
and taking into account that $a(u_n-v_n)\big|_{u_n=v_n-\hbar}=a(-\hbar)=0$ we
conclude that the unique possibility for a non-trivial contribution is:
$\beta_{nn}=\gamma_{nn}=-1$, $\gamma_{nj}=-1$,
$\beta_{nj}=+1$, $j=1,\ldots,n-1$, $\beta_{in}=-1$, $\gamma_{in}=+1$, $i=1,\ldots,n-1$. The last formulae
impose the same DWBC for the $(n-1)\times(n-1)$ sublattice: $\delta_{n-1,j}=\beta_{nj}=+1$, $j=1,\ldots,n-1$,
$\alpha_{i,n-1}=\gamma_{in}=+1$, $i=1,\ldots,n-1$, $d_{n-1,n-1}=d_{nn}$. Thus the substitution $u_n=v_n-\hbar$
to the partition function for the whole lattice yiefds
\begin{eqnarray}
\fl Z^{+-}_{-+}(u_n=v_n-\hbar,u_{n-1},\ldots,u_1;v_n,v_{n-1},\ldots,v_1;\lambda) \ny \\
\fl\quad =\bar c(-\hbar;\lambda)\prod_{j=1}^{n-1}\bar b(v_n-v_j-\hbar;\lambda+(n-j)\hbar)
 \prod_{i=1}^{n-1}b(u_i-v_n;\lambda+(n-i)\hbar) \ny \\
\fl\quad \times Z^{+-}_{-+}(u_{n-1},\ldots,u_1;v_{n-1},\ldots,v_1;\lambda). \label{proofZrec}
\end{eqnarray}
Using the explicit expressions \eqref{BWa} for the Boltzmann weights,
one can rewrite
 the last formula in the form~\eqref{Zrec}. \qed

%{\bfseries Remark \remnumber\label{rem2}.}
\begin{remark}\label{rem2} From Formula~\eqref{proofZrec}
 we see that the following transformation of the $R$-matrix
\begin{equation} \label{remZrec}
 b(u,v;\lambda)\to\rho\,b(u,v;\lambda), \qquad \bar b(u,v;\lambda)\to\rho^{-1}\, \bar b(u,v;\lambda)
\end{equation}
does not change the recursion relation~\eqref{Zrec}, where $\rho$ is a
non-zero constant which does not
depend on $u$, $v$ and $\lambda$.
\end{remark}

\begin{lemma} \label{theor_3prop}
If the set of functions $\{Z^{(n)}(u_n,\ldots,u_1;v_n,\ldots,v_1;\lambda)\}_{n\ge1}$ satisfies
the conditions of Propositions~\ref{propZsym}, \ref{propZep}, \ref{propZrec} and the initial condition
\begin{equation}
 Z^{(1)}(u_1;v_1;\lambda)=\bar c(u_1-v_1)=\frac{\theta(u_1-v_1-\lambda)
 \theta(\hbar)}{\theta(-\lambda)} \label{Z0c}
\end{equation}
then
\begin{equation}
 Z^{+-}_{-+}(u_n,\ldots,u_1;v_n,\ldots,v_1;\lambda)=Z^{(n)}(u_n,\ldots,u_1;v_n,\ldots,v_1;\lambda).
 \label{Z_Z}
\end{equation}
\end{lemma}

\noindent  Due to \eqref{Z0c}, this lemma can be proved by induction
on $n$. Let the equality~\eqref{Z_Z} be valid for $n-1$. Consider the functions $Z^{+-}_{-+}(u_n,
\ldots,u_1;v_n,\ldots,v_1;\lambda)$ and $Z^{(n)}(u_n,\ldots,u_1;v_n,\ldots,v_1;\lambda)$ as functions
of $u_n$. Both are elliptic polynomials of degree $n$ with character~\eqref{char_chi}. They have the same
value at the point $u_n=v_n-\hbar$, and due to the symmetry of these functions with respect to the
parameters $\{v_j\}_{j=1}^n$ they coincide at all points
$u_n=v_j-\hbar$, $j=1,\ldots,n$. It follows from Lemma~\ref{ep_coin}
(see~\ref{Appendix_ep}) that these functions are identical. \qed

%{\bfseries Remark \remnumber\label{rem3}.}
\begin{remark}\label{rem3}
As we can see from the proof of Lemma~\ref{theor_3prop},
it is sufficient to establish the symmetry only with respect to the variables $v_j$.
\end{remark}

%{\bfseries Remark \remnumber\label{rem4}.}
\begin{remark}\label{rem4}
The transformation~\eqref{remZrec} of the $R$-matrix does
not change the partition function with DWBC.
\end{remark}

\section{Elliptic projections of currents}
\label{sec4}

 Let $\lfK_0=\mathbb C[u^{-1}][[u]]$
be the completed space
%?
of complex-valued meromorphic  functions defined in the neighborhood
of the origin which
have only simple poles at this point.  Let
 $\{\epsilon^i\}$ and $\{\epsilon_i\}$ be dual bases in $\lfK_0$ such
 that
 $\oint\frac{du}{2\pi i} \epsilon^i(u)\,\epsilon_j(u)=\delta^i_j$.
%  Below to keep a reasonable size of the paper we are using much material of our
% papers  \cite{PRS} and \cite{EPR05}.

\subsection{Current description of the elliptic algebra}
\label{sec41}

Let $\A$ be a Hopf algebra generated by elements $\hat h[s]$, $\hat e[s]$, $\hat f[s]$,
$s\in\lfK_0$, subject to the linear relations
\begin{equation*}
\hat x[\alpha_1s_1+\alpha_2s_2]=\alpha_1\hat x[s_1]+\alpha_2\hat x[s_2], \qquad \alpha_1,\alpha_2\in\mathbb C,
\quad s_1,s_2\in\lfK_0,
\end{equation*}
where $x\in\{h,e,f\}$. The commutation relations will be written in
terms of the currents,
\begin{eqnarray}
 h^+(u)=\sum_{i\ge0}\hat h[\epsilon^{i;0}]\epsilon_{i;0}(u),\qquad
 h^-(u)=-\sum_{i\ge0}\hat h[\epsilon_{i;0}]\epsilon^{i;0}(u), \ny \\
 f(u)=\sum_i\hat f[\epsilon^i]\epsilon_i(u),\qquad  e(u)=\sum_i\hat e[\epsilon^i]\epsilon_i(u). \label{ftdef}
\end{eqnarray}
The currents $e(u)$ and $f(u)$ are called the {\it total
  currents}. They are defined in terms of dual bases of
$\lfK_0$ and their definition does not depend on the choice
of these dual bases (see \cite{ER1,PRS}).
 The currents $h^+(u)$ and $h^-(u)$ are called the {\it Cartan
   currents} and they are defined in terms of the special basis
\begin{equation*}
 \epsilon^{k;0}(u)=\frac1{k!}\left(\frac{\theta'(u)}{\theta(u)}\right)^{(k)}, \quad k\ge0;
 \qquad % \label{eps0^p} \\
 \epsilon_{k;0}(u)=(-u)^k,\quad k\ge0.   % \label{eps0_p}
\end{equation*}
The commutation relations are~\cite{EF}:
\begin{equation*}
[K^\pm(u),K^\pm(v)]=0, \qquad [K^+(u),K^-(v)]=0,
\end{equation*}
\begin{equation*}
K^\pm(u)e(v)K^\pm(u)^{-1}=\frac{\theta(u-v+\hbar)}{\theta(u-v-\hbar)}e(v),
\end{equation*}
\begin{equation*}
K^\pm(u)f(v)K^\pm(u)^{-1}=\frac{\theta(u-v-\hbar)}{\theta(u-v+\hbar)}f(v), % \label{cr_EFR0KK}
\end{equation*}
\begin{equation}
\theta(u-v-\hbar)e(u)e(v)=\theta(u-v+\hbar)e(v)e(u),\label{cr_EFR0ee}
\end{equation}
\begin{equation}
\theta(u-v+\hbar)f(u)f(v)=\theta(u-v-\hbar)f(v)f(u),\label{cr_EFR0ff}
\end{equation}
\begin{equation*}
[e(u),f(v)]=\hbar^{-1}\delta(u,v)\Big(K^+(u)-K^-(v)\Big),
\end{equation*}
where
$K^+(u)=\exp\Big(\frac{e^{\hbar\partial_u}-e^{-\hbar\partial_u}}{2\partial_u}
h^+(u)\Big)$, $K^-(u)=\exp\big(\hbar h^-(u)\big)$ and $\delta(u,v)=\sum\limits_{n\in\mathbb Z}
\frac{u^n}{v^{n+1}}$ is a delta-function~\footnote{One can find more details about distributions acting on
$\lfK_0$ and their significance in the theory of current algebras in our previous paper~\cite{PRS}.} for $\lfK_0$.
The algebra $\A$ is a non-central version of the algebra $A(\tau)$ introduced in~\cite{ER1}.
This algebra is equipped with the co-product and co-unit:
\begin{equation*}
 \Delta K^\pm(u)=K^\pm(u)\otimes K^\pm(u),
 \end{equation*}
\begin{equation*}
\Delta e(u)=e(u)\otimes1+K^-(u)\otimes e(u),
\end{equation*}
\begin{equation*}
 \Delta f(u)=f(u)\otimes K^+(u)+1\otimes f(u), % \label{cpKdef}
 \end{equation*}
\begin{equation*}
\varepsilon(K^\pm(u))=1, \qquad \varepsilon(e(u))=0, \qquad \varepsilon(f(u))=0.
\end{equation*}

Let $\A_F$ and $\A_E$ be the subalgebras of $\A$ generated by the generators $\hat h[\epsilon_{i;0}]$, $\hat f[s]$,
and $\hat h[\epsilon^{i;0}]$, $\hat e[s]$, respectively, $s\in\lfK_0$. The subalgebra $\A_F$ is
described by the currents $K^+(u)$, $f(u)$, and the subalgebra $\A_E$
by $K^-(u)$, $e(u)$. We
introduce the notation $H^+$ for the subalgebra of $\A$ generated by $\hat h[\epsilon_{i;0}]$.
As stated in~\cite{EF}, the bialgebras $(\A_F,\Delta^{op})$ and
$(\A_E,\Delta)$ are dual
with respect to the Hopf pairing $\langle\cdot,\cdot\ra\colon\A_F\times\A_E\to\mathbb C$ defined
in terms of currents as follows:
\begin{equation}
 \La f(u),e(v)\Ra=\hbar^{-1}\delta(u,v), \qquad
 \La K^+(u),K^-(v)\Ra=\frac{\theta(u-v-\hbar)}{\theta(u-v+\hbar)}. \label{HPdef}
\end{equation}
These formulae uniquely define a Hopf pairing on  $\A_F\times\A_E$. In particular, one can
derive the following formula
\begin{eqnarray}
 \La f(t_n)\cdots f(t_1),e(v_n)\cdots e(v_1)\Ra
  =   \ny \\
  =\hbar^{-n}\sum_{\sigma\in S_n}\prod_{{l<l' \atop \sigma(l)>\sigma(l')}}
   \frac{\theta(v_{\sigma(l)}-v_{\sigma(l')}+\hbar)}{\theta
   (v_{\sigma(l)}-v_{\sigma(l')}-\hbar)}\prod_{m=1}^n\delta(t_m,v_{\sigma(m)}).
    \label{la_ff_ee_ra}
\end{eqnarray}

\subsection{Projections of currents}
\label{sec42}

We define the projections as linear maps acting on the subalgebra
$\A_F$.
Dual projections, which we do not consider here, act in the subalgebra
$\A_E$.
We define the projections in terms of the {\itshape half-currents} $f^+_\lambda(u)$ and
$f^-_\lambda(u)$, defined below.
These are usually defined as parts of the sum~\eqref{ftdef} (with the corresponding sign)
such that $f(u)=f^+_\lambda(u)-f^-_\lambda(u)$. Here $\lambda$ is the
parameter
for
%?
the decomposition
of the total current into the difference of half-currents.
Elliptic half-currents are investigated in details on the classical level in~\cite{PRS}.
We will  introduce the half-currents by their representations by means
of integral
 transforms of the total current $f(u)$:
\begin{equation}
\fl \f^+(u)=\oint\limits_{|v|<|u|}\frac{dv}{2\pi i}\frac{\theta(u-v-\lambda)}
 {\theta(u-v)\theta(-\lambda)}f(v), \quad
 \f^-(u)=\oint\limits_{|v|>|u|}\frac{dv}{2\pi i}\frac{\theta(u-v-\lambda)}
 {\theta(u-v)\theta(-\lambda)}f(v),\label{fpmht}
\end{equation}
where $\lambda\notin \Gamma=\mathbb Z+\mathbb Z\tau$.
The half-current $f^+_\lambda(u)$ is called positive
 and $f^-_\lambda(u)$ is called negative.

The corresponding positive and negative projections are also parameterized by $\lambda$
and they are defined on the half-currents as follows:
\begin{eqnarray}
 P^+_\lambda\big(f^+_\lambda(u)\big)=f^+_\lambda(u),\qquad P^-_\lambda\big(f^+_\lambda(u)\big)=0,
 \label{Pfpdef} \\
 P^+_\lambda\big(f^-_\lambda(u)\big)=0,\qquad  P^-_\lambda\big(f^-_\lambda(u)\big)=f^-_\lambda(u).
 \label{Pfmdef}
\end{eqnarray}

Let us first define the projections in the subalgebra $\A_f$ generated
by the currents $f(u)$.
As a linear space this subalgebra is spanned by the products $f(u_n)f(u_{n-1})\cdots f(u_1)$,
$n=0,1,2,\ldots$. It means that any element of $\A_f$ can be represented as a sum (maybe infinite)
of integrals~\footnote{The integral $\oint$ without limits means a formal integral -- a continuous
extension of the integral over the unit circle.}
\begin{equation*}
 \oint\frac{du_n\cdots du_1}{(2\pi i)^n}f(u_n)\cdots f(u_1)\,s_n(u_n)\cdots s_1(u_1), \qquad
 s_n,\ldots, s_1\in\lfK_0.
\end{equation*}
It follows from the PBW theorem proved in~\cite{EF} that any element
of $\A_f$ can also be represented
as the sum of the integrals
\begin{equation*}
\fl\oint\limits\frac{du_n\cdots du_1}{(2\pi i)^n}
 f^-_{\lambda+2(n-1)\hbar}(u_n)\cdots f^-_{\lambda+2m\hbar}(u_{m+1})
 f^+_{\lambda+2(m-1)\hbar}(u_m)\cdots f^+_{\lambda}(u_1)s_n(u_n)\cdots s_1(u_1),
\end{equation*}
$s_n\cdots s_1\in\lfK_0$, $0\le m\le n$. Therefore, it is sufficient to define the projections
on these products of half-currents:
\begin{equation}
 P^+_\lambda(x^-x^+)=\varepsilon(x^-)x^+,\qquad P^-_\lambda(y^-y^+)=y^-\varepsilon(y^+), \label{Pxy_def}
\end{equation}
where
\begin{eqnarray}
\fl x^-=f^-_{\lambda+2(n-1)\hbar}(u_n)\cdots f^-_{\lambda+2m\hbar}(u_{m+1}),\qquad
y^-=f^-_{\lambda}(u_n)\cdots f^-_{\lambda-2(n-m-1)\hbar}(u_{m+1}),  \ny \\
\fl x^+=f^+_{\lambda+2(m-1)\hbar}(u_m)\cdots f^+_{\lambda}(u_1), \qquad
y^+=f^+_{\lambda-2(n-m)\hbar}(u_m)\cdots f^+_{\lambda-2(n-1)\hbar}(u_1). \ny
\end{eqnarray}
The product of zero number of currents is identified with $1$ and in this case: $\varepsilon(1)=1$.
The counit $\varepsilon$ of a nonzero number of half-currents is always zero. So, this definition
generalizes Formulae~\eqref{Pfpdef} and \eqref{Pfmdef}. We complete the definition of the projections
on the subalgebra $\A_F=\A_f\cdot H^+$ by the formulae
\begin{equation*}
P^+_\lambda(at^+)=P^+_\lambda(a)t^+, \qquad P^-_\lambda(at^+)=P^-_\lambda(a)\varepsilon(t^+),
\end{equation*}
where $a\in\A_f$, $t^+\in H^+$.

\subsection{The projections and the universal elliptic weight function}
\label{sec43}

Consider the expressions of the form
\begin{equation}
 P^+_{\lambda-(n-1)\hbar}\big(f(u_n)f(u_{n-1})\cdots f(u_2)f(u_1)\big), \label{Pfn0}
\end{equation}
where the parameter $\lambda-(n-1)\hbar$ is chosen for symmetry reasons. Let us begin with the
case $n=1$. Formula~\eqref{Pfpdef} implies that in this case the projection is equal to the positive
half-current, which can be represented as an integral transform of the total current:
\begin{equation*}
 P^+_{\lambda}\big(f(u_1)\big)=f^+_{\lambda}(u_1)=
 \oint\limits_{|u_1|>|v_1|}\frac{dv_1}{2\pi i}\frac{\theta(u_1-v_1{-\lambda})}{\theta(u_1-v_1)
 \theta(-\lambda)} f(v_1).
 \end{equation*}
The kernel of this transform gives the initial condition for the partition function with a factor:
\begin{equation}
 Z^{(1)}(u_1;v_1;{\lambda})=
  {\theta(\hbar)\theta(u_1-v_1)}\frac{\theta(u_1-v_1{-\lambda})}{\theta(u_1-v_1)\theta({-\lambda})}.
  \label{icZ0}
\end{equation}

The projections~\eqref{Pfn0} can be calculated by generalizing the method proposed  in~\cite{KhP}
for the
algebra $U_q(\hat\slt)$. The method uses a recursion over $n$. Let us first present the last total
current in~\eqref{Pfn0} as the difference of half-currents:
\begin{eqnarray}
\fl P^+_{\lambda-(n-1)\hbar}\big(f(u_n)\cdots f(u_2)f(u_1)\big)=   \label{Pfn0pm}\\
\fl  =P^+_{\lambda-(n-3)\hbar}\big(f(u_n)\cdots f(u_2)\big)f^+_{\lambda-(n-1)\hbar}(u_1)
   -P^+_{\lambda-(n-1)\hbar}\big(f(u_n)\cdots f(u_2)f^-_{\lambda-(n-1)\hbar}(u_1)\big) . \ny
\end{eqnarray}
In the first term we move out the positive half-current from the projection and, therefore,
calculation of this term reduces to the computation of the $(n-1)$-st projection.
In the second term in~\eqref{Pfn0pm} we move the negative half-current to the left step by step
using the following commutation relation~\cite{EF}
\begin{equation*}
  f(v)f^-_{\lambda}(u_1)
   =\frac{\theta(v-u_1-\hbar)}{\theta(v-u_1+\hbar)}f^-_{\lambda+2\hbar}(u_1)f(v)
     +\frac{\theta(v-u_1+\lambda+\hbar)}{\theta(v-u_1+\hbar)}F_{\lambda}(v), \\
\end{equation*}
where
\begin{equation*}
   F_\lambda(v)=\frac{\theta(\hbar)}{\theta(\lambda+\hbar)}
        \big(f^+_{\lambda+2\hbar}(v)f^+_{\lambda}(v)
                               -f^-_{\lambda+2\hbar}(v)f^-_{\lambda}(v)\big).
\end{equation*}

At each step we obtain an additional term containing $F_\lambda(u)$ and at the last step the negative
half-current is annihilated by the projection:
\begin{equation}
   P^+_{\lambda-(n-1)\hbar}\big(f(u_n)\cdots f(u_2)f^-_{\lambda-(n-1)\hbar}(u_1)\big)
       =\sum_{j=2}^n Q_j(u_1) X_j,  \label{Pfn0m}
\end{equation}
where
\begin{eqnarray*}
   Q_j(u)=\frac{\theta(u_j-u+\lambda-(n-2j+2)\hbar)}{\theta(u_j-u+\hbar)}
        \prod_{k=2}^{j-1}\frac{\theta(u_k-u-\hbar)}{\theta(u_k-u+\hbar)}, \\
   X_j=P^+_{\lambda-(n-1)\hbar}\big(f(u_n)\cdots f(u_{j+1})
            F_{\lambda-(n-2j+3)\hbar}(u_j)f(u_{j-1})\cdots f(u_2)\big).
\end{eqnarray*}
Setting $u_1=u_i$ in~\eqref{Pfn0m}, we can substitute the negative
half-current for the positive one
using the commutation relation for the total currents $f(u)$ and the equality $f(u)f(u)=0$.
Moving out the positive half-current to the left one obtains a linear system of equations for
$X_i$,  $i=2,\ldots,n$:
\begin{equation}
   P^+_{\lambda-(n-3)\hbar}\big(f(u_n)\cdots f(u_2)\big)f^+_{\lambda-(n-1)\hbar}(u_i)
                =\sum_{j=2}^n Q_j(u_i)X_j. \label{PfQX}
\end{equation}
Multiplying each equation~\eqref{PfQX} by
\begin{equation*}
 \frac{\theta(u_i-u+\lambda)}{\theta(\lambda)}
  \prod_{k=2}^n\frac{\theta(u_k-u_i+\hbar)}{\theta(u_k-u+\hbar)}
  \prod_{{k=2 \atop k\ne i}}^n\frac{\theta(u_k-u)}{\theta(u_k-u_i)},
\end{equation*}
summing over $i=2,\ldots,n$ and using the interpolation formula (see~\ref{Appendix_ep})
\begin{equation}
 Q_j(u)=\sum_{i=2}^n Q_j(u_i)\frac{\theta(u_i-u+\lambda)}{\theta(\lambda)}
  \prod_{k=2}^n\frac{\theta(u_k-u_i+\hbar)}{\theta(u_k-u+\hbar)}
  \prod_{{k=2 \atop k\ne i}}^n\frac{\theta(u_k-u)}{\theta(u_k-u_i)} \label{interpol_Q_p}
\end{equation}
yields
\begin{eqnarray}
\fl  P^+_{\lambda-(n-3)\hbar}\big(f(u_n)\cdots f(u_2)\big)
   \sum_{i=2}^n \frac{\theta(u_i-u+\lambda)}{\theta(\lambda)}
  \prod_{k=2}^n\frac{\theta(u_k-u_i+\hbar)}{\theta(u_k-u+\hbar)}
  \prod_{{k=2 \atop k\ne i}}^n\frac{\theta(u_k-u)}{\theta(u_k-u_i)}f^+_{\lambda-(n-1)\hbar}(u_i)
 \ny\\
    =\sum_{j=2}^n Q_j(u)X_j. \label{Pn1fQX}
\end{eqnarray}
Comparing~\eqref{Pn1fQX} with~\eqref{Pfn0m}, we conclude that
\begin{eqnarray}
\fl   P^+_{\lambda-(n-1)\hbar}\big(f(u_n)\cdots f(u_2)f^-_{\lambda-(n-1)\hbar}(u_1)\big)
  =P^+_{\lambda-(n-3)\hbar}\big(f(u_n)\cdots f(u_2)\big) \ny \\
\fl\quad   \times\sum_{i=2}^n\frac{\theta(u_i-u_1+\lambda)}{\theta(\lambda)}
  \prod_{k=2}^n\frac{\theta(u_k-u_i+\hbar)}{\theta(u_k-u_1+\hbar)}
  \prod_{{k=2 \atop k\ne i}}^n\frac{\theta(u_k-u_1)}{\theta(u_k-u_i)}f^+_{\lambda-(n-1)\hbar}(u_i).
\end{eqnarray}
Finally, returning to Formula~\eqref{Pfn0pm} we derive the following expression for
the projection
\begin{equation}
\fl   P^+_{\lambda-(n-1)\hbar}\big(f(u_n)\cdots f(u_2)f(u_1)\big)
  =P^+_{\lambda-(n-3)\hbar}\big(f(u_n)\cdots f(u_2)\big)f^+_{\lambda-(n-1)\hbar}(u_1;u_n,\ldots,u_2),
  \label{Pfn0_Pfn1fp}
\end{equation}
where we introduce the linear combination of the  currents:
\begin{eqnarray}
\fl f^+_{\lambda-(n-2m+1)\hbar}(u_m;u_n,\ldots,u_{m+1})= f^+_{\lambda-(n-2m+1)\hbar}(u_m)
 -\sum_{i=m+1}^n\frac{\theta(u_i-u_m+\lambda+(m-1)\hbar)}{\theta(\lambda+(m-1)\hbar)} \ny \\
\fl \times\prod_{k=m+1}^n \frac{\theta(u_k-u_i+\hbar)}{\theta(u_k-u_m+\hbar)}
  \prod\limits_{{k=m+1 \atop k\ne i}}^n\frac{\theta(u_k-u_m)}{\theta(u_k-u_i)}
   f^+_{\lambda-(n-2m+1)\hbar}(u_i). \label{fpm_nm}
\end{eqnarray}
Continuing this computation by induction we obtain an expression for the projections in terms of the
half-currents~\eqref{fpm_nm}:
\begin{equation}
 \fl  P^+_{\lambda-(n-1)\hbar}\big(f(u_n)\cdots f(u_2)f(u_1)\big)
   =\prod_{n\ge m\ge1}^{\longleftarrow}f^+_{\lambda-(n-2m+1)\hbar}(u_m;u_n,\ldots,u_{m+1}). \label{Pfn0_prod_f}
\end{equation}
Using the addition formula
\begin{equation}
\prod_{i=1}^n G_{\lambda_i}(u_i-v)=\sum_{i=1}^n\prod_{{j=1 \atop j\ne i}}^n
              G_{\lambda_j}(u_j-u_i)G_{\lambda_0}(u_i-v), \label{GGG}
\end{equation}
where $G_{\lambda}(u-v)=\frac{\theta(u-v+\lambda)}{\theta(u-v)\theta(\lambda)}$,
$\lambda_0=\sum\limits_{i=1}^n\lambda_i$, one can represent the half-currents~\eqref{fpm_nm}
as integral transforms of the total current:
\begin{eqnarray}
\fl f^+_{\lambda-(n-2m+1)\hbar}(u_m;u_n,\ldots,u_{m+1})= \ny \\
\fl= \prod_{k=m+1}^n\frac{\theta(u_k-u_m)}{\theta(u_k-u_m+\hbar)}
  \oint\limits_{|u_i|>|v|}\frac{dv}{2\pi i}
   \frac{\theta(u_m-v-\lambda-(m-1)\hbar)}
                      {\theta(u_m-v)\theta(-\lambda-(m-1)\hbar)}
     \prod_{k=m+1}^n\frac{\theta(u_k-v+\hbar)}{\theta(u_k-v)}
                                            f(v). \ny
\end{eqnarray}
Replacing each combination of the half-currents~\eqref{fpm_nm}  in~\eqref{Pfn0_prod_f}
by their integral form we obtain
\begin{eqnarray}
\fl P^+_{\lambda-(n-1)\hbar}\big(f(u_n)\cdots f(u_2)f(u_1)\big)
 =\prod_{n\ge k>m\ge1}\frac{\theta(u_k-u_m)}{\theta(u_k-u_m+\hbar)}
   \oint\limits_{|u_i|>|v_j|}\frac{dv_n\cdots dv_1}{(2\pi i)^n} \ny \\
\fl     \prod_{n\ge k>m\ge1}\frac{\theta(u_k-v_m+\hbar)}{\theta(u_k-v_m)}
    \prod_{m=1}^n \frac{\theta(u_m-v_m-\lambda-(m-1)\hbar)}{\theta(u_m-v_m)\theta(-\lambda-(m-1)\hbar)}
           f(v_n)\cdots f(v_1). \label{Pfn0_prod_f_int}
\end{eqnarray}

Formulae \eqref{Pfn0_prod_f} and  \eqref{Pfn0_prod_f_int}
yield expressions for the universal elliptic weight
functions in terms of the current generators of the algebra $\mathcal{A}$.

\subsection{Universal weight function and SOS model partition function}

To extract the kernel from the expression \eqref{Pfn0_prod_f_int} and derive a formula for the partition function
we use the Hopf pairing~\eqref{HPdef}. Let us calculate the following expression generalizing~\eqref{icZ0}:
\begin{eqnarray}
\fl Z^{(n)}(u_n,\ldots,u_1;v_n,\ldots,v_1;\lambda)=\prod_{i,j=1}^n\theta(u_i-v_j)
  \prod_{n\ge k>m\ge1}\frac{\theta(u_k-u_m+\hbar)\theta(v_k-v_m-\hbar)}{\theta(u_k-u_m)
  \theta(v_k-v_m)}\ny \\
\fl\quad \times\big(\hbar\theta(\hbar)\big)^n\La P^+_{\lambda-(n-1)\hbar}\big(f(u_n)\cdots f(u_1)\big),
e(v_n)\cdots e(v_1)\Ra. \label{Zn_def}
\end{eqnarray}
Using the expression for the projection of the product of the total currents~\eqref{Pfn0_prod_f_int}
and Formula~\eqref{la_ff_ee_ra} we obtain
\begin{eqnarray}
\fl Z^{(n)}(u_n,\ldots,u_1;v_n,\ldots,v_1;\lambda)= \ny \\
\fl\quad   =\theta(\hbar)^n\prod_{i,j=1}^n\theta(u_i-v_j)
   \prod_{k>m}\frac{\theta(v_k-v_m-\hbar)}{\theta(v_k-v_m)}
  \sum_{\sigma\in S_n}\prod_{{l<l' \atop \sigma(l)>\sigma(l')}}
   \frac{\theta(v_{\sigma(l)}-v_{\sigma(l')}+\hbar)}{\theta(v_{\sigma(l)}-v_{\sigma(l')}-\hbar)}\times \ny \\
\fl\qquad  \times\prod_{k>m}\frac{\theta(u_k-v_{\sigma(m)}+\hbar)}{\theta(u_k-v_{\sigma(m)})}
   \prod_{m=1}^n\frac{\theta(u_m-v_{\sigma(m)}-\lambda-(m-1)\hbar)}
                      {\theta(u_m-v_{\sigma(m)})\theta(-\lambda-(m-1)\hbar)}= \label{Zn_theta} \\
\fl\quad =\prod_{k>m}\frac{\theta(v_k-v_m-\hbar)}{\theta(v_k-v_m)}
  \sum_{\sigma\in S_n}\prod_{{l<l' \atop \sigma(l)>\sigma(l')}}
   \frac{\theta(v_{\sigma(l)}-v_{\sigma(l')}+\hbar)}{\theta(v_{\sigma(l)}-v_{\sigma(l')}-\hbar)}\times \ny \\
\fl\qquad  \times\prod_{k>m}\theta(u_k-v_{\sigma(m)}+\hbar)\prod_{k<m}\theta(u_k-v_{\sigma(m)})
   \prod_{m=1}^n\frac{\theta(u_m-v_{\sigma(m)}-\lambda-(m-1)\hbar)\theta(\hbar)}
                      {\theta(-\lambda-(m-1)\hbar)}. \ny
\end{eqnarray}
From this formula we see that the expression~\eqref{Zn_theta} defines a holomorphic function of
the variables $u_i$.

\begin{theorem}
The set of functions $\big\{Z^{(n)}(u_n,\ldots,u_1;v_n,\ldots,v_1;\lambda)\big\}_{n\ge1}$
defined by Formula~\eqref{Zn_def} satisfies the conditions of
Propositions~\ref{propZsym}, \ref{propZep}, \ref{propZrec} and the initial condition~\eqref{Z0c}.
%Therefore, by virtue of the theorem~\ref{theor_3prop}
They coincide with the partition functions of
the SOS model with DWBC:
\begin{eqnarray}
\fl Z^{+-}_{-+}(u_n,\ldots,u_1;v_n,\ldots,v_1;\lambda)= \ny \\
\fl\quad =\prod_{n\ge k>m\ge1}\frac{\theta(v_k-v_m-\hbar)}{\theta(v_k-v_m)}
  \sum_{\sigma\in S_n}\prod_{{l<l' \atop \sigma(l)>\sigma(l')}}
   \frac{\theta(v_{\sigma(l)}-v_{\sigma(l')}+\hbar)}{\theta(v_{\sigma(l)}-v_{\sigma(l')}-\hbar)}
     \prod_{1\le k<m\le n}\theta(u_k-v_{\sigma(m)}) \ny \\
\fl\quad\times\prod_{n\ge k>m\ge1}\theta(u_k-v_{\sigma(m)}+\hbar)\prod_{m=1}^n
\frac{\theta(u_m-v_{\sigma(m)}-\lambda-(m-1)\hbar)\theta(\hbar)}
                      {\theta(-\lambda-(m-1)\hbar)}. \label{Z_Z_}
\end{eqnarray}
\end{theorem}

\noindent
The initial condition~\eqref{Z0c} is satisfied because
Formula~\eqref{icZ0} is satisfied.
%?
The first factor in the right-hand side of~\eqref{Zn_def} is symmetric with respect to both sets
of variables. The symmetry with respect to the variables $\{u\}$ and the variables $\{v\}$ follows
from the commutation relations~\eqref{cr_EFR0ff} and \eqref{cr_EFR0ee} respectively.
Formula~\eqref{Zn_theta} implies that~\eqref{Zn_def} are elliptic polynomials of degree $n$
with character~\eqref{char_chi} in the variables $u_i$, in particular
in $u_n$. We now substitute
$u_n=v_n-\hbar$ to~\eqref{Zn_theta}. The non-vanishing terms in the right-hand side correspond
to the permutations $\sigma\in S_n$ satisfying $\sigma(n)=n$. Substituting $u_n=v_n-\hbar$ into these
 terms, one obtains the recursion relation~\eqref{Zrec}. \qed \\

\section{Degeneration of the partition function}
\label{sec5}

In this section, we investigate the trigonometric degenerations of the formulae obtained
in the elliptic case. In particular, taking the corresponding trigonometric limit in the
expression for the SOS model partition function~\eqref{Z_Z_} reproduces the expression for the
6-vertex partition function~\eqref{stat-s_pr}.

First we consider the degeneration of the $R$-matrix, the matrix of Boltzmann weights,
which defines the model. To do so we need the formula for the trigonometric degeneration
($\tau\to i\infty$) of the odd theta function,
%? defined by the conditions~\r{theta}:
\begin{equation*}
 \lim_{\tau\to i\infty}\theta(u)=\frac{\sin\pi u}{\pi}.
\end{equation*}
In terms of the multiplicative variables $z=e^{2\pi iu}$, $w=e^{2\pi iv}$, this formula can be
rewritten as follows:
\begin{equation*}
2\pi i e^{\pi i(u+v)}\lim_{\tau\to i\infty}\theta(u-v)=z-w.
\end{equation*}
Multiplying the $R$-matrix~\eqref{Rz} by $2\pi i e^{\pi i(u+v)}$ and taking the limit we
obtain the following matrix which depends rationally on the multiplicative variables $z$, $w$
and on the multiplicative parameters $q=e^{\pi i\hbar}$, $\mu=e^{2\pi i\lambda}$:
\begin{eqnarray}
 R(z,w;\mu)=2\pi i e^{\pi i(u+v)}\lim_{\tau\to i\infty}R(u-v;\lambda)= \ny \\
 =\left(\begin{array}{cccc}
  zq-wq^{-1}  & 0              & 0              & 0              \\
  0 & \frac{(z-w)(\mu q-q^{-1})}{(\mu-1)} &\frac{(z-w\mu)(q-q^{-1})}{(1-\mu)}  & 0              \\
  0 & \frac{(z\mu-w)(q-q^{-1})}{(\mu-1)} &\frac{(z-w)(\mu q^{-1}-q)}{(\mu-1)}  & 0   \\
  0              & 0              & 0              & (zq-wq^{-1})
 \end{array}\right).  \label{Rzmu}
\end{eqnarray}
The matrix~\eqref{Rzmu} inherits the property of satisfying the dynamical Yang-Baxter equation
and it defines a statistical model called the {\it trigonometric SOS model}.

To obtain the non-dynamical trigonometric case we need to implement the additional
limit $\lambda\to-i\infty$ implying $\mu\to\infty$ (or $\lambda\to i\infty$ implying $\mu\to0$):
\begin{equation}
\fl \tilde R(z,w)=\lim_{\mu\to\infty}R(z,w;\mu)=\left(\begin{array}{cccc}
  zq-wq^{-1}  & 0              & 0              & 0              \\
  0 & q(z-w) &(q-q^{-1})w  & 0              \\
  0 & (q-q^{-1})z & q^{-1}(z-w)  & 0   \\
  0              & 0              & 0              & zq-wq^{-1}
 \end{array}\right). \label{Rzmund}
\end{equation}
The matrix~\eqref{Rzmund} differs from the matrix of the Boltzmann
weights of the 6-vertex model~\eqref{Bw1}
%?
by the transformation~\eqref{remZrec}. Taking into account
Remark~\ref{rem4} (from Subsection 3.2),
we conclude that both matrices~\eqref{Bw1} and \eqref{Rzmund} define
the same partition function $Z(\{z\},\{w\})$ with DWBC~\footnote{The matrix~\eqref{Rzmund} is the
limit of a matrix which differs from~\eqref{Rz} by the transformation~\eqref{remZrec} with $\rho=q$.}.

To obtain the partition function with DWBC for the trigonometric SOS model,
 one should multiply the partition function with DWBC for the
elliptic SOS model by a factor  and take the trigonometric limit:
\begin{eqnarray}
\fl Z^{+-}_{-+}(\{z\},\{w\};\mu)
 =\prod_{k,j=1}^n\big(2\pi i e^{\pi i(u_k+v_j)}\big)\lim_{\tau\to i\infty}
 Z^{+-}_{-+}(\{u\},\{v\};\lambda)= \ny \\
\fl\quad  =\prod_{n\ge k>m\ge1}\frac{w_k q^{-1}-w_m q}{w_k-w_m}
  \sum_{\sigma\in S_n}\prod_{{l<l' \atop \sigma(l)>\sigma(l')}}
   \frac{w_{\sigma(l)}q-w_{\sigma(l')}q^{-1}}{w_{\sigma(l)}q^{-1}-w_{\sigma(l')}q} \label{Z_Z_t} \\
\fl  \times\prod_{n\ge k>m\ge1}\big(z_k q-w_{\sigma(m)}q^{-1}\big)
\prod_{1\le k<m\le n}\big(z_k-w_{\sigma(m)}\big)
   \prod_{m=1}^n\frac{\big(z_m-w_{\sigma(m)}\mu q^{2(m-1)}\big)(q-q^{-1})}
                      {\big(1-\mu q^{2(m-1)}\big)}. \ny
\end{eqnarray}
It is easy to prove that Formula~\eqref{stat-s_pr} is obtained from Formula~\eqref{Z_Z_t}
by taking the limit: $Z(\{z\},\{w\})=\lim\limits_{\mu\to\infty} Z^{+-}_{-+}(\{z\},\{w\};\mu)$.

\section*{Acknowledgements}
This paper is part of the PhD thesis
of A.S. who he is preparing it under the co-direction of S. P. and V. R. in
the Bogoliubov  Laboratory of Theoretical Physics, JINR, Dubna and in LAREMA,
D\'epartement de Math\'ematiques, Universit\'e d'Angers. He is
grateful to the CNRS-Russia exchange program on mathematical physics  and  personally to J.-M. Maillet
for financial and general support of this thesis project. V. R. is
grateful to
Ph.~Di~Franchesco and T.~Miwa for their stimulating lectures and their interest during
the ENIGMA School on ``Quantum Integrability'' held in Lalonde-les-Maures on October 14-19, 2007.
He thanks O.~Babelon and M.~Talon for the invitation to this school.
We would like to thank Y.~Kosmann-Schwarzbach for her attention, help and useful remarks.
During this project V. R. used
partial financial support from ANR GIMP
and support from the INFN-RFBR ``Einstein grant'' (Italy-Russia).
S. P. was supported in part by RFBR grant 06-02-17383. Both V. R. and
S. P. were supported in part
by the grant for support of scientific schools NSh-8065.2006.2.

\appendix

\section{Interpolation formula for elliptic polynomials}
\label{Appendix_ep}

A group homomorphism $\chi\colon\Gamma\to\mathbb C^\times$,
where $\Gamma=\mathbb Z+\tau\mathbb Z$ and $\mathbb C^\times$ is the multiplicative group of nonzero
complex numbers is called a {\it character}. Each character $\chi$ and integer number
$n$ define a space $\Theta_n(\chi)$ consisting of the holomorphic functions on $\mathbb C$
with the translation properties
\begin{equation*}
 \phi(u+1)=\chi(1)\phi(u), \qquad
 \phi(u+\tau)=\chi(\tau)e^{-2\pi inu-\pi in\tau}\phi(u).
\end{equation*}
If $n>0$ then $\dim\Theta_n(\chi)=n$ (and $\dim\Theta_n(\chi)=0$ if $n<0$). The elements
of the space $\Theta_n(\chi)$ are called elliptic polynomials (or theta-functions) of degree $n$
with character $\chi$.

\begin{prop}
Let $\{\phi_j\}_{j=1}^n$ be a basis of $\Theta_n(\chi)$, with character
$\chi(1)=(-1)^n$, $\chi(\tau)\hm=(-1)^n e^{2\pi i \alpha}$, then the determinant of
the matrix $||\phi_j(u_i)||_{\leq i,j \leq n}$ is equal to
\begin{equation}
  \det||\phi_j(u_i)||=C\cdot\theta(\sum_{k=1}^n u_k-\alpha)\prod_{i<j}\theta(u_i-u_j), \label{ell_Vand}
\end{equation}
where $C$ is a nonzero constant.
\end{prop}
\noindent
Consider the ratio
\begin{equation}
  \frac{\det||\phi_j(u_i)||}{\theta(\sum_{k=1}^n u_k-\alpha)\prod_{i<j}\theta(u_i-u_j)}. \label{ell_Vand_ratio}
\end{equation}
This is an elliptic function of each $u_i$ with only simple poles in any fundamental domain
(the points $u_i$ satisfying $\sum_{k=1}^n u_k-\alpha\in\Gamma$).
Therefore, it is a  constant function of each $u_i$.
Thus this ratio does not
depend on $u_i$ and we have to prove that it does not vanish, that is
that the determinant $\det||\phi_j(u_i)||$ is not identically zero. Let us denote by
 $\Delta^{i_1,\ldots,i_k}_{j_1,\ldots,j_k}$ the minor of this determinant corresponding to
 the $i_1$-th, $\ldots$, $i_k$-th rows and the $j_1$-th, $\ldots$, $j_k$-th columns. Suppose
 that this determinant is identically zero and consider the following decomposition
\begin{equation}
 \det||\phi_j(u_i)||=\sum_{k=1}^n(-1)^{k+1}\phi_k(y_1)\Delta^{2,\ldots,n}_{1,\ldots,k-1,k+1,\ldots,n}.
 \label{ell_Vand_decomp}
\end{equation}
Since the functions $\phi_k(y_1)$ are linearly independent, the minors
$\Delta^{2,\ldots,n}_{1,\ldots,k-1,k+1,\ldots,n}$ are identically zero. Decomposing
the minor $\Delta^{2,\ldots,n}_{2,\ldots,n}$ we conclude that the minors
 $\Delta^{3,\ldots,n}_{2,\ldots,k-1,k+1,\ldots,n}$ are identically zero, and so on.
 Finally, we obtain that $\Delta^{n}_{n}=\phi_n(y_n)$ is identically zero which cannot be true. \qed

\begin{lemma} \label{ep_coin}
 Let us consider two elliptic polynomials $P_1,P_2\in\Theta_n(\chi)$,
where $\chi(1)=(-1)^n$, $\chi(\tau)\hm=(-1)^n e^{\alpha}$, and $n$ points
$u_i$, $i=1,\ldots,n$, such that $u_i-u_j\not\in\Gamma$, $i\ne j$,
and $\sum_{k=1}^n u_k-\alpha\not\in\Gamma$. If the values of these polynomials
coincide at these points, $P_1(u_i)=P_2(u_i)$, then these polynomials coincide:
$P_1(u)=P_2(u)$.
\end{lemma}

\noindent  Decomposing the polynomials under consideration as
$P_a(u)=\sum_{i=1}^n p_a^i\phi_i(u)$, $a=1,2$, we obtain the system of equations
\begin{equation*}
 \sum_{i=1}^n p_{12}^i\phi_i(u)=0,
\end{equation*}
with respect to the variables $p_{12}^i=p_{1}^i-p_2^i$. We have proved
that
the determinant of this system is equal to~\eqref{ell_Vand}
and therefore is not zero. Hence, this system has only the trivial solution $p_{12}^i=0$,
but this implies $P_1(u)=P_2(u)$. \qed

Let $P\in\Theta_n(\chi)$ be an elliptic polynomial, where $\chi(1)=(-1)^n$,
$\chi(\tau)\hm=(-1)^n e^{2\pi i\alpha}$, and $u_i$, $i=1,\ldots,n$, be $n$ points
such that $u_i-u_j\not\in\Gamma$, $i\ne j$,
and $\sum_{k=1}^n u_k-\alpha\not\in\Gamma$. This polynomial can be
recovered from the values at these points:
\begin{equation}
\fl\qquad P(u)=\sum_{i=1}^n P(u_i)
   \frac{\theta(u_i-u+\alpha-\sum_{m=1}^n u_m)}{\theta(\alpha-\sum_{m=1}^n u_m)}
  \prod_{{k=1 \atop k\ne i}}^n\frac{\theta(u_k-u)}{\theta(u_k-u_i)}. \label{interp_P}
\end{equation}
Indeed, the right hand side belongs to $\Theta_n(\chi)$, this equality
holds at the points $u=u_i$.
Using Lemma~\ref{ep_coin}, we conclude that~\eqref{interp_P} holds at all $u\in\mathbb C$.

Consider the meromorphic functions
\begin{equation*}
 \fl\qquad  Q_j(u)=\frac{\theta(u_j-u+\lambda-(n-2j+2)\hbar)}{\theta(u_j-u+\hbar)}
        \prod_{k=2}^{j-1}\frac{\theta(u_k-u-\hbar)}{\theta(u_k-u+\hbar)}.
\end{equation*}
It is easy to check that the functions
\begin{eqnarray*}
 \fl\quad  P_j(u)=\prod_{k=2}^n\theta(u_k-u+\hbar)Q_j(u)= \ny\\
  \fl\qquad       =\theta(u_j-u+\lambda-(n-2j+2)\hbar)\prod_{k=2}^{j-1}\theta(u_k-u-\hbar)
           \prod_{k=j+1}^n\theta(u_k-u+\hbar)
\end{eqnarray*}
belong to $\Theta_{n-1}(\chi)$, where $\chi(1)=(-1)^{n-1}$, $\chi(\tau)\hm=(-1)^{n-1}
e^{2\pi i\alpha}$, $\alpha=\lambda+\sum_{k=2}^n u_k$. Since $\lambda\not\in\Gamma$,
the polynomials $P_j(u)$ can be recovered from by their values $P_j(u_i)$ via the interpolation
 formula~\eqref{interp_P}. Taking into account the relation between $Q_j(u)$ and
 $P_j(u)$ we obtain Formula~\eqref{interpol_Q_p}.~\footnote{We
 can require the condition $u_i-u_j\not\in\Gamma$,
 because the $u_i$'s in Formula~\eqref{interpol_Q_p} are formal variables.}

\subsection*{References}
\bibliographystyle{amsplain}

\end{document}